# VARIATIONS AND EXTENSIONS OF CROFT'S PROBLEM

SOPHIA SMYRLI

ABSTRACT. In this work we study the following classical still challenging Calculus problem: If $f : (0, \infty) \to \mathbb{R}$ is a continuous function, for which the sequence $\{f(nx)\}$ tends to zero, for every positive $x$, as $n$ tends to infinity, then $f(x)$ also tends to zero, as $x$ tends to infinity.

## 1. INTRODUCTION

**PROBLEM.** *Let $f : (0, \infty) \to \mathbb{R}$ be a continuous function, and suppose that $\lim_{n \to \infty} f(nx) = 0$, for all $x > 0$. Then $\lim_{x \to \infty} f(x) = 0$.*

This problem first appears in 1957, in a paper by H. T. Croft [1]. In this work we provide variations and extensions of this problem. Some of these have appeared in the paper of J.F.C. Kingsman [3].

Croft's proof is based on the following lemma[1]:

**Lemma 1.** *If $U$ is an open and unbounded subset of $(0, \infty)$, then there exists an $x > 0$, such that $nx \in U$, for infinitely many $n \in \mathbb{N}$.*

Using this lemma, if $f(x)$ does not tend to zero, as $x \to \infty$, then there is an $\varepsilon > 0$, such that at least one of the sets

$$U = \{x \in (0, \infty) : f(x) > \varepsilon\}, \quad V = \{x \in (0, \infty) : f(x) < -\varepsilon\}$$

is unbounded. Without loss of generality assume that $U$ is unbounded, and it is also open, Lemma 1 applies, and hence there exists an $x \in (0, \infty)$, such that $nx \in U$, for infinitely many $n \in \mathbb{N}$. This means that the sequence $f(nx)$ does not tends to zero.

*Second proof based on Baire's Theorem.* Let $\varepsilon > 0$ arbitrary. We define the sets

$$F_n = \{x \in (0, \infty) : |f(mx)| \leq \varepsilon, \quad \text{for all } m \geq n\}, \quad n \in \mathbb{N} = \{1, 2, \ldots\},$$

which are closed in $(0, \infty)$. Since $f(nx) \to 0$, for every $x > 0$, then

$$\bigcup_{n \in \mathbb{N}} F_n = (0, \infty).$$

Baire's Theorem (see Lemma 2 and Remark 3 in Section 3) provides that at least one the $F_n$'s has nonempty interior. In particular, there exist $n \in \mathbb{N}$ and $a, b > 0$, with $a < b$, such that

$$[a, b] \subset F_n.$$

This means that, for every $m \geq n$,

$$x \in [ma, mb] \quad \implies \quad |f(x)| \leq \varepsilon.$$

Next we observe that if $k > a/(b-a)$ then

$$kb > (k+1)a$$

---

and thus
$$[ka, kb] \cup [(k+1)a, (k+1)b] \cup [(k+2)a, (k+2)b] \cup \cdots = [ka, \infty).$$
Consequently, if also $k \geq n$, then $[ka, \infty) \subset f^{-1}\big([-\varepsilon, \varepsilon]\big)$, and hence
$$x \geq ka \implies |f(x)| \leq \varepsilon,$$
which finally implies that $\lim_{x \to \infty} f(x) \to 0$. $\square$

**Remark 1.** *Counterexamples*
  (i) *In general, this result does not hold if the function $f : (0, \infty) \to \mathbb{R}$ is not continuous. For example, let*
$$f(x) = \begin{cases} 1 & \text{if } x = p, \text{ where } p \text{ is a prime,} \\ 0 & \text{otherwise.} \end{cases}$$
  *Then every $x > 0$, the sequence $\{f(nx)\}$, eventually vanishes, and hence*
$$\lim_{n \to \infty} f(nx) = 0,$$
  *while $f(p) = 1$, and hence $f(x)$ does not tend to zero, as $x$ tends to infinity.*
  (ii) *In the next example we see that if $f$ is continuous and $f(nx) \to 0$, for all $x$ in a dense subset of $(0, \infty)$, then $f(x)$ does not necessarily tend to zero, as $x \to \infty$.*
  *Let $f : (0, \infty) \to \mathbb{R}$, be defined as*
$$f(x) = \sin(n!\pi x), \quad \text{whenever} \quad x \in [n-1, n].$$
  *Clearly, if $x = p/q$, with $p, q \in \mathbb{N}$, then*
$$f(nx) = 0, \quad \text{for all} \quad n \geq q,$$
  *and hence $\lim_{n \to \infty} f(nx) = 0$, whenever $x$ is rational. Meanwhile, $\max_{x \in [n-1, n]} f(x) = 1$, and hence the limit $\lim_{x \to \infty} f(x)$ does not exist.*

## 2. Variations and extensions

  (i) Let $f : (0, \infty) \to \mathbb{R}$, be a continuous function, for which the limit $\lim_{n \to \infty} f(nx)$ exists for all $x > 0$, then the limit $\lim_{x \to \infty} f(x)$, also exists.

*Proof.* Our target is to prove that both $\liminf_{x \to \infty} f(x)$ and $\limsup_{x \to \infty} f(x)$ are finite, and their difference
$$\limsup_{x \to \infty} f(x) - \liminf_{x \to \infty} f(x)$$
is arbitrarily small. Let $\varepsilon > 0$ arbitrary. We define the sets
$$F_{m,j} = \big\{x \in (0, \infty) : f(nx) \in [(j-1)\varepsilon, (j+1)\varepsilon], \text{ for all } n \geq m\big\}, \quad m \in \mathbb{N}, \ j \in \mathbb{Z}, \qquad (2.1)$$
which are closed in $(0, \infty)$. We have that,
$$\bigcup_{j \in \mathbb{Z}} \bigcup_{m \in \mathbb{N}} F_{m,j} = (0, \infty).$$
To see this, let $x \in (0, \infty)$. As $f(nx) \to a$, for some $a \in \mathbb{R}$, there exists an $m \in \mathbb{N}$ such that
$$n \geq m \implies f(nx) \in (a - \varepsilon/2, a + \varepsilon/2) \subset [(j-1)\varepsilon, (j+1)\varepsilon],$$
where $j = \left[\frac{a}{\varepsilon} + \frac{1}{2}\right]$. Baire's Theorem provides that at least one the $F_{m,j}$'s has nonempty interior. In particular, there exist $m \in \mathbb{N}$, $j \in \mathbb{Z}$, and $a, b > 0$, with $a < b$, such that $[a, b] \subset F_{m,j}$. As in the case of the original problem
$$k > \max\left\{m, \frac{a}{b-a}\right\} \implies [ka, \infty) \subset f^{-1}\big([(j-1)\varepsilon, (j+1)\varepsilon]\big),$$

which implies that
$$(j-1)\varepsilon \le \liminf_{x\to\infty} f(x) \le \limsup_{x\to\infty} f(x) \le (j+1)\varepsilon.$$

Thus $\liminf_{x\to\infty} f(x), \limsup_{x\to\infty} f(x) \in \mathbb{R}$, and $\limsup_{x\to\infty} f(x) - \liminf_{x\to\infty} f(x) \le 2\varepsilon$.

(*ii*) In fact, we do not need to restrict the limits of the sequences $\{f(nx)\}$ in $\mathbb{R}$:

Let $f : (0, \infty) \to \mathbb{R}$, be a continuous function, such that for all $x > 0$, the limit $\lim_{n\to\infty} f(nx)$ exists in $[-\infty, \infty]$, then the limit $\lim_{x\to\infty} f(x)$, also exists $[-\infty, \infty]$.

*Proof.* Here $[-\infty, \infty] = \mathbb{R} \cup \{-\infty, \infty\}$. We reduce the problem to the previous one by defining
$$g(x) = \frac{f(x)}{1+|f(x)|}.$$
Then $g(x) \in (-1, 1)$, and
$$f(x) = \frac{g(x)}{1-|g(x)|}.$$
Clearly, the sequence $\{g(nx)\}$, converges for all $x > 0$, and and its limit lies in $[-1, 1]$. In particular, if $\lim_{n\to\infty} g(nx) = \ell$, then
$$\lim_{n\to\infty} f(nx) = \begin{cases} \dfrac{\ell}{1-|\ell|} & \text{if } |\ell| < 1, \\ \infty & \text{if } \ell = 1, \\ -\infty & \text{if } \ell = -1. \end{cases}$$

The previous example provides that the limit $\lim_{x\to\infty} g(x)$ exists, and in fact lies in $[-1, 1]$. Consequently, the limit $\lim_{x\to\infty} f(x)$ exists in $[-\infty, \infty]$.

(*iii*) There exists a variation for inferior and superior limits.

If $f : (0, \infty) \to \mathbb{R}$ is continuous and $\limsup_{n\to\infty} f(nx) = s$, (*resp.* $\limsup_{n\to\infty} f(nx) \le s$), for all $x \in (0, \infty)$, then $\limsup_{x\to\infty} f(x) = s$ (*resp.* $\limsup_{x\to\infty} f(x) \le s$).

First note that $\limsup_{n\to\infty} f(nx) = s$, means that for every $\varepsilon > 0$,
(a.) Eventually $f(nx) < s + \varepsilon$, and
(b.) $f(nx) > s - \varepsilon$, for infinitely many $n \in \mathbb{N}$.
As before, we define
$$F_n = \{x \in (0, \infty) : f(mx) < s + \varepsilon, \text{ for all } m \ge n\}.$$

Clearly, $\bigcup_{n\in\mathbb{N}} F_n = (0, \infty)$, and hence $F_n$, for some $n \in \mathbb{N}$, contains an non-trivial interval, say $[a, b]$. This means that $[ma, mb] \subset F_n$, for all $m \ge n$, and thus
$$[ka, \infty) \subset f^{-1}(-\infty, s + \varepsilon),$$
for $k > \max\{n, a/(b-a)\}$, which means that $\limsup_{x\to\infty} f(x) < s + \varepsilon$. Meanwhile, the fact that $\limsup_{n\to\infty} f(n) = s$, implies that $\limsup_{x\to\infty} f(x) \ge s$.

(*iv*) If the limit of $f(x)$, as $x$ tends to infinity, does not exist, then as $f$ is continuous, for every
$$\liminf_{x\to\infty} f(x) \le s \le \limsup_{x\to\infty} f(x),$$
there exists an increasing sequence $\{x_n\} \subset (0, \infty)$, with $x_n \to \infty$ and $f(x_n) \to s$. In particular, we can pick a specific type of such sequence:

*If $f : (0, \infty) \to \mathbb{R}$ is continuous and $\liminf_{x \to \infty} f(x) \le s \le \limsup_{x \to \infty} f(x)$, then there exists an $x \in (0, \infty)$ and an increasing sequence $\{k_n\}$ of positive integers, such that $\lim_{n \to \infty} f(k_n x) = s$.*

(v) A further extension is to replace the range of $f$ by an arbitrary metric space:

*Let $(X, d)$ be a metric space and $f : (0, \infty) \to X$. If the limit $\lim_{n \to \infty} f(nx)$ exists for all $x > 0$, then the limit $\lim_{x \to \infty} f(x)$ also exists.*

Let $\hat{X}$ be the completion of $X$. For simplicity, we assume that $X \subset \hat{X}$ and the metric $d$ is extended to $\hat{X}$ with the same symbol. Observe that the image of $f$:
$$Y = f(0, \infty) = \bigcup_{n \in \mathbb{N}} f[2^{-n}, n] \subset X,$$
is $\sigma$–compact, and hence separable. Let $S$ be a countable and dense subset of $Y$ and $\varepsilon > 0$. We define
$$F_{n,s} = \{x \in (0, \infty) : d(f(mx), s) \le \varepsilon, \text{ for all } m \ge n\}, \quad n \in \mathbb{N}, \ s \in S.$$
Using the argument of the original problem we obtain that $[a, b] \subset F_{n,s}$, for some $n, j \in \mathbb{N}$, and hence, for every $\varepsilon > 0$, there exist an $s = s_\varepsilon \in S$ and $M = M_\varepsilon > 0$, such that
$$x \ge M_\varepsilon \implies f(x) \in \overline{B}(s_\varepsilon, \varepsilon) = \{s \in \hat{X} : d(s_\varepsilon, s) \le \varepsilon\}.$$
Clearly, if $x \ge \max\{M_{1/k} : k = 1, \ldots, n\}$, then
$$f(x) \in K_n = \bigcap_{j=1}^{n} \overline{B}(s_{1/j}, 1/j).$$
Clearly, $f(x) \in K_n$, for all $x \ge R_n = \max\{M_{1/k} : k = 1, \ldots, n\}$, and hence $K_n \ne \emptyset$. Also
$$\text{diam}(K_n) \le \text{diam}(\overline{B}(s_{1/n}, 1/n)) \le \frac{2}{n}.$$
Completeness of $\hat{X}$ implies that $\bigcap_{n \in \mathbb{N}} K_n \ne \emptyset$. In fact, $\bigcap_{n \in \mathbb{N}} K_n$ is a singleton. Let $s_0$ its unique element. We have that
$$x \ge R_n \implies d(f(x), s_0) \le \frac{2}{n}.$$
Thus $\lim_{x \to \infty} f(x) = s_0$ and $s_0 = \lim_{n \to \infty} f(n) \in X$.

(vi) The range of $f$ can become even more abstract. Namely, our result holds even when the range of $f$ is a regular Hausdorff space.

**Definition 1.** A topological space $X$ is called *regular* if given any closed subset $F$ of $X$ and any point $x \in X \setminus F$, there exist an open neighborhood $U$ of $x$ and an open neighborhood $V$ of $F$, which are disjoint. A topological space is called regular Hausdorff or simply $T_3$ if it is both Hausdorff and regular.

We have the following extension:

*Let $X$ be a regular Hausdorff topological space, and $f : (0, \infty) \to X$ be a continuous map. If the limit $\lim_{n \to \infty} f(nx)$ exists for all $x > 0$, then the limit $\lim_{x \to \infty} f(x)$ also exists.*

*Proof.* Let's remind that if $\{t_n\}$ is a sequence in $X$, this sequence converges to $t \in X$, if for every open neighborhood $U$ of $t$, there exists an $n_0 \in \mathbb{N}$, such that
$$n \ge n_0 \implies t_n \in U.$$
The assumption that $X$ is a Hausdorff space guarantees the uniqueness of the limit.

At first, we shall show that all the sequences $\{f(nx)\}$ converge to the same limit. If this were not the case, then there would exist $x_1, x_2 > 0$, such that

$$t_1 = \lim_{n \to \infty} f(nx_1) \neq \lim_{n \to \infty} f(nx_2) = t_2.$$

Since $X$ is a Hausdorff space, we can pick open and disjoint neighborhoods $V_1$ and $V_2$, of $t_1$ and $t_2$, respectively. Let

$$U_1 = f^{-1}[V_1] \quad \text{and} \quad U_2 = f^{-1}[V_2].$$

The sets $U_1$ and $U_2$ are open, since $f$ is continuous. Also, $U_1$ and $U_2$ are unbounded, since $f(nx_i) \in V_i$, for $i = 1, 2$, and all $n$ sufficiently large. Using Lemma 4, we can obtain an $x \in (0, \infty)$ and a strictly increasing sequence $\{k_n\} \subset \mathbb{N}$, such that

$$f(k_{2n-1}x) \in V_1 \quad \text{while} \quad f(k_{2n}x) \in V_2.$$

We shall construct $\{k_n\}$ and a *decreasing* sequence of nontrivial closed intervals $[a_n, b_n]$, recursively. Let $[a_0, b_0] \subset (0, \infty)$ be an arbitrary nontrivial interval. Assuming that we have already defined the terms $k_1, \ldots, k_{n-1}$, and the intervals $[a_1, b_1], \ldots, [a_{n-1}, b_{n-1}]$, then using Lemma 4, we shall define $k_n$ and $[a_n, b_n]$ as follows. If $n$ is odd, then we pick a term $k_n > k_{n-1}$ and a nontrivial interval $[a_n, b_n]$, with $[a_{n-1}, b_{n-1}] \subset [a_n, b_n]$, so that

$$k_n[a_n, b_n] \subset V_1,$$

whereas if $n$ is even, then the $k_n$ and the $[a_n, b_n]$ are picked so that

$$k_n[a_n, b_n] \subset V_2.$$

Clearly, $\bigcap_{n \in \mathbb{N}} [a_n, b_n] \neq \emptyset$. Let $x \in \bigcap_{n \in \mathbb{N}} [a_n, b_n]$. Then, for every $n \in \mathbb{N}$,

$$\ell_{k_{2n-1}} x \in V_1 \quad \text{while} \quad \ell_{k_{2n}} x \in V_2,$$

and hence the limit $\lim_{n \to \infty} f(nx)$ does not exist. Therefore, all the sequences $\{f(nx)\}$, $x \in (0, \infty)$, converge to the same limit $t \in X$.

It now remains to show that $\lim_{x \to \infty} f(x) = t$. In particular, we need to show that for every open neighborhood $V$ of $t$, there exists an $M > 0$, such that

$$x > M \implies f(x) \in V.$$

Since $X$ is regular Hausdorff, $F = X \setminus V$ is closed and $x \notin F$, there exist $U$ and $W$ disjoint open sets such that $x \in U$ and $F \subset W$. In fact,

$$U \subset X \setminus W \subset V \implies U \subset \overline{U} \subset \overline{X \setminus W} = X \setminus W \subset V,$$

since $X \setminus W$ is closed. We set

$$F_n = \{x \in (0, \infty) : f(mx) \in \overline{U}, \quad \text{for all } m \geq n\},$$

which are closed in $(0, \infty)$. Since $f(nx) \to t$, for every $x > 0$, then

$$\bigcup_{n \in \mathbb{N}} F_n = (0, \infty).$$

Baire's Theorem provides that at least one the $F_n$'s has nonempty interior. In particular, there exist $n \in \mathbb{N}$ and $a, b > 0$, with $a < b$, such that

$$[a, b] \subset F_n.$$

As we have seen before, this implies that,

$$x \geq ka \implies f(x) \in \overline{U} \subset V,$$

provided that $k > a/(b - a)$, and finally implies that $\lim f(x) \to t$. □

(*vii*) Another variation is to replace the domain of $f$ by $(1, \infty)$, and to replace the given converging sequences by sequences of the form $\{f(x^n)\}$, for all $x \in (1, \infty)$.

(*viii*) Another variation is to replace the domain of $f$ by $\mathbb{R}$, and to replace the given converging sequences by sequences of the form $\{f(x + \log n)\}$, for all $x \in \mathbb{R}$.

(*ix*) Another variation is to replace the domain of $f$ by $\mathbb{C} \smallsetminus \{0\}$. For example we could have the following result:

Let $f : \mathbb{C} \smallsetminus \{0\} \to \mathbb{C}$ be a continuous function and assume that the limit $\lim_{n \to \infty} f(nz)$, exists, for all $z \in \mathbb{C} \smallsetminus \{0\}$. Then the limit $g(\vartheta) = \lim_{t \to \infty} f(te^{i\vartheta})$ also exists, for all $\vartheta \in [0, 2\pi]$, and $g$ is continuous in a dense $G_\delta$ subset of $[0, 2\pi]$.

This is due to the fact that $g$ is a pointwise limit of continuous functions.

(*x*) Seeking plausible potential generalizations:

**Question.** *Let $\{\ell_n\}$ an increasing sequence in $(0, \infty)$, with $\ell_n \to \infty$ and $f : (0, \infty) \to \mathbb{R}$ be a continuous function, such that the limit $\lim_{n \to \infty} f(\ell_n x)$ exists for all $x \in (0, \infty)$. Does this imply that the limit $\lim_{x \to \infty} f(x)$ also exists?*

*Answer.* Not is general. For example if $\ell_n = a^n$, for $a > 1$, and
$$f(x) = \sin(2\pi \log_a x), \quad x \in (0, \infty),$$
then
$$f(\ell_n x) = f(a^n x) = \sin\left(2\pi \log_a(a^n x)\right) = \sin\left(2\pi \log_a x + 2\pi n\right) = \sin\left(2\pi \log_a x\right) = f(x),$$
for all $n \in \mathbb{N}$. A counterexample can be constructed of a function for which $f(a^n x) \to 0$, for all $x$, but the limit $\lim_{x \to \infty} f(x)$ does not exist.

Nevertheless, the answer becomes positive in $\ell_n$ grows sub-exponentially:

**Improved result.** *Let $\{\ell_n\}$ be a strictly increasing sequence in $(0, \infty)$, with $\ell_n \to \infty$ and*
$$\lim_{n \to \infty} \frac{\ell_{n+1}}{\ell_n} = 1,$$
*and let $f : (0, \infty) \to \mathbb{R}$ be a continuous function, such that the limit $\lim_{n \to \infty} f(\ell_n x)$ exists for all $x \in (0, \infty)$. Then the limit $\lim_{x \to \infty} f(x)$ also exists.*

*Proof.* Initiating the proof of the first extension of the original problem, we prove that the difference $\limsup_{x \to \infty} f(x) - \liminf_{x \to \infty} f(x)$ is arbitrarily small. In order to achieve this we define the following sets
$$F_{n,j} = \left\{x \in (0, \infty) : f(\ell_m x) \in [(j-1)\varepsilon, (j+1)\varepsilon], \text{ for all } m \geq n\right\}, \quad n \in \mathbb{N}, \ j \in \mathbb{Z},$$
which are closed in $(0, \infty)$. Once again $\bigcup_{j \in \mathbb{Z}} \bigcup_{n \in \mathbb{N}} F_{n,j} = (0, \infty)$. Baire's Theorem provides the existence of $n \in \mathbb{N}$, $j \in \mathbb{Z}$, and $0 < a < b$, such that $[a, b] \subset F_{n,j}$. This means that
$$\bigcup_{m \geq n} [ma, mb] \subset f^{-1}\left([(j-1)\varepsilon, (j+1)\varepsilon]\right).$$
Meanwhile, if $k > \max\left\{n, \dfrac{a}{b-a}\right\}$, then $\bigcup_{m \geq k} [ma, mb] = [ka, \infty)$, and hence
$$x \geq ka \implies f(x) \in [(j-1)\varepsilon, (j+1)\varepsilon],$$
and therefore
$$\limsup_{x \to \infty} f(x) - \liminf_{x \to \infty} f(x) \leq 2\varepsilon.$$

*Second proof.* We shall use Lemma 4. If $\lim_{x \to \infty} f(x)$ does not exist, then we can pick $A, B \in \mathbb{R}$, with
$$\liminf_{x \to \infty} f(x) < A < B < \limsup_{x \to \infty} f(x),$$
and set
$$V = \{x \in (0, \infty) : f(x) < A\} \quad \text{and} \quad W = \{x \in (0, \infty) : f(x) > B\}.$$
Since $f$ is continuous, $V$ and $W$ are both open and also unbounded. Next, we shall define a subsequence $\{\ell_{k_n}\}$ of $\{\ell_n\}$ and a sequence of nontrivial intervals $[a_n, b_n] \subset (0, \infty)$, $n \in \mathbb{N}$, such that $[a_{n+1}, b_{n+1}] \subset [a_n, b_n]$,
$$\ell_{k_{2n-1}}[a_{2n-1}, b_{2n-1}] \subset V \quad \text{and} \quad \ell_{k_{2n}}[a_{2n}, b_{2n}] \subset W, \quad \text{for all } n \in \mathbb{N}.$$
This can be achieved recursively, starting with an arbitrary nontrivial interval $[a_0, b_0] \subset (0, \infty)$. Assuming that we have already defined the terms $\ell_{k_1}, \ldots, \ell_{k_n}$, and the intervals
$$[a_1, b_1], \ldots, [a_{n-1}, b_{n-1}],$$
then using Lemma 4, we shall define $\ell_{k_n}$ and $[a_n, b_n]$ as follows. If $n$ is odd, then we pick a term $\ell_{k_n} > \ell_{k_{n-1}}$ and a nontrivial interval $[a_n, b_n]$, with $[a_{n-1}, b_{n-1}] \subset [a_n, b_n]$, and
$$\ell_{k_n}[a_n, b_n] \subset V,$$
whereas if $n$ is even, then the $\ell_{k_n}$ and the $[a_n, b_n]$ are picked so that
$$\ell_{k_n}[a_n, b_n] \subset W.$$
Clearly, $\bigcap_{n \in \mathbb{N}}[a_n, b_n] \neq \emptyset$. Let $x \in \bigcap_{n \in \mathbb{N}}[a_n, b_n]$. Then, for every $n \in \mathbb{N}$,
$$\ell_{k_{2n-1}}x \in V \quad \text{while} \quad \ell_{k_{2n}}x \in W.$$
Therefore
$$f(\ell_{k_{2n-1}}x) < A < B < f(\ell_{k_{2n}}x), \quad \text{for all } n \in \mathbb{N},$$
and hence the limit $\lim_{n \to \infty} f(\ell_n x)$ does not exist. □

($xi$) The previous result implies that:

*If $\{a_n\}$ is an increasing sequence, such that $a_n \to \infty$, while $a_{n+1} - a_n \to 0$, and $f : \mathbb{R} \to \mathbb{R}$ is continuous, and the limit $\lim_{n \to \infty} f(x + a_n)$, exists for every $x \in \mathbb{R}$, then the limit $\lim_{x \to \infty} f(x)$ also exists.*

**Remark 2.** *Kingsman [3] and Fehér et al [2] considered the characterization of sequences $\{c_n\}$, with the property that if $f : \mathbb{R} \to \mathbb{R}$ is continuous and $\lim_{n \to \infty} f(x + c_n)$ exists for all $x \in \mathbb{R}$, then the limit $\lim_{x \to \infty} f(x)$ also exists.*

## 3. Our tools

**Lemma 2. (Baire Category Theorem)** *Let $(X, d)$ be a complete metric space and $V_n$, $n \in \mathbb{N}$, a sequence of open and dense subsets of $X$. Then their intersection $\bigcap_{n \in \mathbb{N}} V_n$ is also dense in $X$.*

*Proof.* See [5, page 11].

An equivalent formulation of Baire's Theorem is:

**Corollary 1.** *Let $(X, d)$ be a complete metric space and $F_n$, $n \in \mathbb{N}$, a sequence of closed subsets of $X$, such that $\bigcup_{n \in \mathbb{N}} F_n = X$. Then at least one of the $F_n$'s has nonempty interior.*

**Remark 3.** *Both Lemma 2 and Corollary 1 still hold even when the complete metric space $X$ is replaced by an open subset of $X$.*

**Lemma 3.** *Let $U$ be an open and unbounded subset of $(0, \infty)$ and $a, b > 0$, with $a < b$. Then for every $n \in \mathbb{N}$, there exist $n' \in \mathbb{N}$, with $n' \geq n$ and $a', b' > 0$, with $a' < b'$, such that $[a', b'] \subset (a, b)$, and*
$$n'[a', b'] \subset U.$$

Instead of proving this lemma, we shall prove the following slightly stronger result:

**Lemma 4.** *Let $\{\ell_n\} \subset (0, \infty)$, be a strictly increasing sequence which tends to infinity, such that $\lim_{n \to \infty} \ell_{n+1}/\ell_n = 1$. Let also $U$ be an open and unbounded subset of $(0, \infty)$ and $a, b > 0$, with $a < b$. Then for every $n \in \mathbb{N}$, there exist an $n' \geq n$ and $a', b'$, with $a' < b'$ and $[a', b'] \subset (a, b)$, such that*
$$\ell_{n'}[a', b'] \subset U.$$

*Proof.* Since $\lim_{n \to \infty} \ell_{n+1}/\ell_n = 1$, there exists an $n_0$, such that
$$\frac{\ell_{n+1}}{\ell_n} < \frac{b}{a}, \quad \text{for } n \geq n_0.$$
Since $U$ is not bounded, there exists a $t \in U$, such that $t/b > \ell_{n_0}$, and as $\{\ell_n\}$ is strictly increasing and $\ell_n \to \infty$, there exists $n_1 \in \mathbb{N}$, such that
$$\ell_{n_1-1} \leq \frac{t}{b} < \ell_{n_1},$$
Clearly, $n_1 > n_0$, and hence
$$\ell_{n_1} = \ell_{n_1-1} \cdot \frac{\ell_{n_1}}{\ell_{n_1-1}} < \frac{t}{b} \cdot \frac{b}{a} = \frac{t}{a}.$$
Thus
$$\frac{t}{b} < \ell_{n_1} < \frac{t}{a} \quad \implies \quad \ell_{n_1} a < t < \ell_{n_1} b.$$
This implies that the set $U \cap (\ell_{n_1} a, \ell_{n_1} b)$ is non-empty and it is also open. Hence there exists a nontrivial interval
$$[r, s] \subset U \cap (\ell_{n_1} a, \ell_{n_1} b).$$
If we set $a' = r/\ell_{n_1}$ and $b' = s/\ell_{n_1}$, the interval $[\ell_{n_1} a', \ell_{n_1} b'] = [r, s] \subset (a, b)$, and hence $[a', b']$ satisfies our requirements. $\square$

Corresponding author: Sophia Smyrli

Department of Mathematics, University of Cyprus
*E-mail address*: sophia.smyrli@gmail.com

Department of Mathematics & Statistics,, University of Cyprus